\documentclass{amsart}
\input amssym.def
\input amssym
\usepackage[body={7in,9.5in},top=0.5in, left=0.8in ]{geometry}
\usepackage{hyperref}
\newtheorem{theorem}{Theorem}[section]
\newtheorem{definition}[theorem]{Definition}
\newtheorem{lemma}[theorem]{Lemma}
\newtheorem{proposition}[theorem]{Proposition}
\newtheorem{corollary}[theorem]{Corollary}

\newtheorem{remark}[theorem]{Remark}

\theoremstyle{definition}

\theoremstyle{remark}

\numberwithin{equation}{section}

\numberwithin{equation}{section}

\title[General composition of formal power seres -- geometry and topology]{Topology and geometry of the general composition of formal power series -- towards Fr\'{e}chet-Lie group-like formalism}

\author{Dawid Bugajewski\\ F\lowercase{aculty of }P\lowercase{hysics, }U\lowercase{niversity of }W\lowercase{arsaw, }W\lowercase{arsaw, }P\lowercase{oland}}



\address[Dawid Bugajewski]{Faculty of Physics\\
University of Warsaw\\
Pasteura 5\\
02-093 Warsaw\\
Poland}
		
\email[Dawid~Bugajewski]{d.bugajewski@student.uw.edu.pl}

\begin{document}
\maketitle
\begin{abstract}
In this article, we study the properties of the autonomous superposition operator on the space of formal power series, including those with nonzero constant term. We prove its continuity and smoothness with respect to the topology of pointwise convergence and a natural Fr\'{e}chet manifold structure. A necessary and sufficient condition for the left composition inverse of a formal power series to exist is provided. We also present some properties of the Fr\'{e}chet-Lie group structures on the set of nonunit formal power series.
\end{abstract}

\section{Introduction}
One of the most important operations on the space of formal power series is their composition, first defined under the assumption that the inner series is a nonunit ($a_0=0$), then without this restriction, see for example the monograph \cite{g21}. The first case seems to be a much simpler one, since no infinite sums of series coefficients have to be considered. It is well known that the set of nonunit formal power series with composition, as well as a narrower set of formal power series with $a_0=0$, $a_1=1$, usually denoted as $\xi(\mathbb C)$ (or, more generally, $\xi(\mathbf{k})$, where $\mathbf{k}$ is a commutative ring), forms a group. The topological and geometrical properties of $\xi(\mathbf{k})$ have been widely investigated for over fifty years, starting from the pioneering papers of Got\^o \cite{goto} or Jennings \cite{jennings}, which, unfortunately, were not initially widely noticed. A comprehensive discussion of this theory can be found e.g. in \cite{babenko}. It is worth mentioning that not only the groups $\xi(\mathbf{k})$, where $\mathbf{k}$ is a field of characteristic zero like $\mathbb C$, but also groups of formal power series over $\mathbf{k}=\mathbb Z$ or $\mathbb Z_p$ are an important subject of investigation, see e.g. \cite{babenko2}, \cite{york}, \cite{camina}.\\

A significant step towards understanding the second case -- the general composition of formal power series, where both series can have nonzero constant terms -- was the necessary and sufficient condition for $g\circ f$ ($g,f\in\mathbb X(\mathbb C)$) to exist provided by Gan and Knox in 2002 \cite{gank}. This result, known as the General Composition Theorem, has become a crucial tool in investigating the properties of general composition, for instance the general chain rule \cite{chain}, the boundary behavior of power series \cite{gb1}, the J.C.P Miller formula \cite{jcp} or a more general case, where the outer series is a formal Laurent series \cite{laurent}. In this article, it will allow us to take steps towards a more topological and geometrical description of the composition of formal power series without imposing any nonunitness conditions. \\

The structure of this paper is the following. First, in Section \ref{2}, we provide some definitions and theorems used in the main part of the article. Then, in Section \ref{3}, we analyze the general autonomous superposition operator on the space of formal power series $\mathbb X(\mathbb C)$ (including those with $a_0\neq 0$). We prove its continuity with respect to the natural pointwise convergence topology and its smoothness as a mapping between two Fr\'{e}chet manifolds. We also provide a necessary and sufficient condition for a left composition inverse of a given formal power series $f=a_0+a_1z+...$ to exist. Finally, we give a short description and present some algebraic properties of the Fr\'{e}chet-Lie group of nonunit formal power series, which we denote $\mathbb X^0_z(\mathbb C)$, a natural extension of the group $\xi(\mathbb C)$ of formal power series satisfying $a_0=0$, $a_1=1$.

\section{Preliminaries}\label{2}

In this section we collect the most important definitions and facts which will be used in the sequel.\\
We will denote by $\mathbb N_0$ (by $\mathbb N$) the set of all nonegative (positive) integers.

\begin{definition}\label{def1}
A formal power series on $\mathbb C$ is defined as a mapping $f:\mathbb N_0\rightarrow \mathbb C$, usually denoted by
\begin{eqnarray*}
f=a_0+a_1z+a_2z^2+\cdots\mbox{ or }f=\sum\limits_{n=0}^{\infty}a_nz^n.
\end{eqnarray*}
If $a_0=0$, then $f$ is called a nonunit, otherwise $f$ is called a unit. We define deg$(f)=\max\{n:a_n\neq 0\}\in\mathbb N_0\cup\{+\infty\}$, ord$(f)=\min\{n:a_n\neq 0\}\in\mathbb N_0\cup\{+\infty\}$, and denote by $\mathbb X(\mathbb C)$ ($\mathbb X^0(\mathbb C)$) the set of all (nonunit) formal power series over $\mathbb C$. The product of two formal power series is defined analogously to the standard Cauchy product, that is
\[
\left(\sum\limits_{n=0}^{\infty}a_nz^n\right)\left(\sum\limits_{n=0}^{\infty}b_nz^n\right)=\sum\limits_{n=0}^{\infty}\left(\sum\limits_{k=0}^na_kb_{n-k}\right)z^n.
\]
\end{definition}
Also, let us denote $[z^n]f:=a_n$ ($n\in\mathbb N_0$). \\Now, it is well known that the multiplicative inverse $f^{-1}$ of a given $f\in\mathbb X(\mathbb C)$ exists if and only if its zeroth coefficient $a_0$ is nonzero (see e.g. \cite{g21}, Th. 1.1.8).
In this paper, however, we are mainly going to deal with another operation -- the {\it composition} of formal power series, which is defined as follows \cite{gank}:
\begin{definition}\label{defcomp}
Let $g=\sum\limits_{k=0}^\infty b_kz^k \in \mathbb X(\mathbb C)$. We define
\[
\mathbb X_g = \left\{f=a_0+a_1z+... \in \mathbb X(\mathbb C):\sum_{k=0}^\infty b_k a_n^{(k)}\in\mathbb C\mbox{ for every }n\in\mathbb N_0\right\}\subset \mathbb X(\mathbb C),
\]
where $a_n^{(k)}$ denote the coefficients of powers $f^k$ of $f$ for $k\in\mathbb N_0$ ($f^0:=1$, $f^1:=f$, $f^{k}:=f\cdot f^{k-1}$ for $k>1$). The mapping $T_g:\mathbb X_g \rightarrow \mathbb X(\mathbb C)$ such that
\[
T_g(f) = \sum_{k=0}^\infty c_k z^k,\mbox{ where }c_k= \sum\limits_{k=0}^\infty b_k a_n^{(k)}
\]
 is well-defined (e.g. $\mathbb X^0(\mathbb C)\subset \mathbb X_g$); we call $g\circ f:=T_g(f)$ the general composition of formal power series $g$ and $f$.
\end{definition}
The necessary and sufficient condition for the existence of the general composition of two formal power series is given in the following
\begin{theorem}\label{gct}\cite{gank}
(The General Composition Theorem) Let $f=a_0+a_1z+a_2z^2+...$ and $g=b_0+b_1z+b_2z^2+...$ be two formal power series over $\mathbb C$ and deg$(f)\neq 0$. Then the composition $g\circ f$ exists, if and only if
\[
\sum\limits_{n=k}^{\infty}\binom{n}{k}b_na_0^{n-k}\in\mathbb C\mbox{ for every }k\in\mathbb N_0,
\]
or, equivalently $g^{(k)}(a_0)\in\mathbb C$ for every $k\in\mathbb N_0$ ($g^{(k)}$ denotes the $k$th order formal derivative of $g$, $g^{(k)}=\sum\limits_{n=k}^{\infty}\frac{n!b_n}{(n-k)!}z^{n-k}$, and $g^{(k)}(a_0):=\sum\limits_{n=k}^{\infty}\frac{n!}{(n-k)!}b_{n}a_0^{n-k}$, i.e. evaluation of $g^{(k)}$ at $a_0$ as if it was a "classical" (not formal) power series).
\end{theorem}
We will denote by $r(g)$ the radius of convergence of a power series $g$. In particular, if $|a_0|<r(g)$, then $g\circ f$ exists \cite{gank}.\\
It should also be noted that the general composition of formal power series is closely related to the issue of boundary convergence of power series -- for more details see \cite{gb1}; here we provide one fact which will be particularly useful in the sequel:
\begin{theorem}\label{l9.2.5}(\cite{gb1}, Lemma 2.5) 
Let $g\in\mathbb X(\mathbb C)$, $0<r(g)<+\infty$. If $g^{(k)}(z_0)\in\mathbb C$ for every $k\in\mathbb N_0$ for some $z_0\in\mathbb C$, $|z_0|=r(g)$, then for every $k\in\mathbb N_0$, the power series $g^{(k)}(z)$ is uniformly continuous on the closed disc $D=\left\{z\in\mathbb C:|z|\leq r(g)\right\}$.
\end{theorem}
Finally, let us provide some basic information about the so called matrix representation of general composition, which will be useful in the analysis of its reciprocity in the next section; for more information, see e.g. \cite{g21}, par. 5.2.
\begin{definition}(\cite{g21}, Def. 5.2.1) 
For every $f=a_0+a_1z+...\in\mathbb X(\mathbb C)$, we define the composition matrix of $f$ by the formula 
\[
C_f:=\left[\begin{array}{llllll}
a_0^{(0)} & a_1^{(0)} & a_2^{(0)} & \ldots & a_n^{(0)} & \ldots  \\
a_0^{(1)} & a_1^{(1)} & a_2^{(1)} & \ldots & a_n^{(1)} & \ldots \\
a_0^{(2)} & a_1^{(2)} & a_2^{(2)} & \ldots & a_n^{(2)} & \ldots \\
\vdots & \vdots & \vdots & \ddots & \vdots & \ddots \\
a_0^{(k)} & a_1^{(k)} & \ldots & \ldots & a_n^{(k)} & \ldots \\
\vdots & \vdots & \ddots & \ddots & \vdots & \ddots \\
\end{array}\right].
\]
\end{definition}
\begin{proposition}\label{matrixcomp}(\cite{g21}, Prop. 5.2.2) 
Let $f=\sum\limits_{n=0}^{\infty}a_nz^n,g=\sum\limits_{n=0}^{\infty}b_nz^n\in\mathbb X(\mathbb C)$. Then $g\circ f=\sum\limits_{n=0}^{\infty}c_nz^n$, where
\[
\left[\begin{array}{l}
c_0\\
c_1\\
c_2\\
\vdots\\
c_k\\
\vdots\\
\end{array}\right]=C_f^T\left[\begin{array}{l}
b_0\\
b_1\\
b_2\\
\vdots\\
b_k\\
\vdots\\
\end{array}\right],
\]
provided that the product on the right-hand side exists. Otherwise $f\notin\mathbb X_g$.
\end{proposition}
\section{Analytical properties of the autonomous superposition operator and the composition inverse in $\mathbb X(\mathbb C)$}\label{3}

At the beginning of this section let us describe the topology with which the set $\mathbb X(\mathbb C)$ will be endowed. We define a mapping $d:\mathbb X^2(\mathbb C)\mapsto[0,+\infty)$ such that for every $f(z)=\sum\limits_{n=0}^{\infty}a_nz^n$, $g(z)=\sum\limits_{n=0}^{\infty}b_nz^n$,
\[
d(f,g)=\sum\limits_{n\in\mathbb N_0}\frac{1}{2^n}\frac{|a_n-b_n|}{|a_n-b_n|+1}.
\]
It is well known (see e.g. \cite{al}) that $(\mathbb X(\mathbb C),d)$ is a complete, separable metric space. Moreover, let $(f_k)_{k\in\mathbb N}$, $f_k=\sum\limits_{n=0}^{\infty}a_{n,k}z^n$ be a sequence of formal power series and let $f(z)=\sum\limits_{n=0}^{\infty}a_nz^n\in\mathbb X(\mathbb C)$. Then $f_k(z)\stackrel{k\rightarrow\infty}{\rightarrow}f$, if and only if $a_{n,k}\stackrel{k\rightarrow\infty}{\rightarrow}a_n$ for every nonnegative integer $n$.\\

Our first result will concern the continuity of the autonomous superposition operator $T_g:\mathbb X_g\mapsto\mathbb X(\mathbb C)$. Before presenting it, let us clarify that $\mathbb X_g$ is not necessarily an open subset of $\mathbb X(\mathbb C)$ -- it is however, by Theorems \ref{gct} and \ref{l9.2.5}, always of the form $\left\{a_0+a_1z+...:a_0\in D\right\}$, where $D$ is either $\left\{0\right\}$, an open disk centered at 0, a closed disk centered at 0 or the whole complex plane. We endow it with the subspace topology, and since $\mathbb X(\mathbb C)$ is a metric space, we therefore can examine $T_g$'s continuity by considering only sequences $(f_n)$ of elements of $\mathbb X_g$ converging to some limit $f\in\mathbb X_g$ and checking the convergence of $(T_g(f_n))$ in $\mathbb X(\mathbb C)$.

\begin{theorem}\label{thm1}
Let $g\in\mathbb X(\mathbb C)$ and let $T_g:\mathbb X_g\ni f\mapsto g\circ f\in\mathbb X(\mathbb C)$. Then
\begin{enumerate}
\item[a)] if $r(g)=0$, then $T_g$ is continuous;
\item[b)] if $r(g)>0$, then $T_g$ is continuous in every point $f=\sum\limits_{n=0}^{\infty}a_nz^n$ such that (1) $|a_0|<r(g)$ or (2) $r(g)<+\infty$, $|a_0|=r(g)$ and $\mbox{deg}(f)>0$.
\end{enumerate}
\end{theorem}

\begin{proof}
Denote $g(z)=\sum\limits_{n=0}^{\infty}b_nz^n$, let $f(z)=\sum\limits_{n=0}^{\infty}a_nz^n\in\mathbb X_g$ and denote $g\circ f(z)=\sum\limits_{n=0}^{\infty}c_nz^n$. Let $(f_k)_{k\in\mathbb N}$, $f_k(z)=\sum\limits_{n=0}^{\infty}a_{n,k}z^n$ be a sequence of elements of $\mathbb X_g$ converging to $f$. We will divide the proof into two cases.
\begin{enumerate}
\item[(a)] $r(g)=0$:\\
Then all elements of $\mathbb X_g$ are nonunits and $g\circ f_k(z)=\sum\limits_{n=0}^{\infty}c_{n,k}z^n$, where
\begin{eqnarray*}
c_{n,k}=\sum\limits_{s=0}^nb_sa_{n,k}^{(s)}=\sum\limits_{s=0}^n\sum\limits_{R_{n,s}}\frac{b_ss!a_{0,k}^{r_0}\ldots a_{n,k}^{r_n}}{r_0!\ldots r_n!}\stackrel{k\rightarrow\infty}{\rightarrow}\sum\limits_{s=0}^n\sum\limits_{R_{n,s}}\frac{b_ss!a_{0}^{r_0}\ldots a_{n}^{r_n}}{r_0!\ldots r_n!}=\sum\limits_{s=0}^nb_sa_{n}^{(s)}=c_n,
\end{eqnarray*}
for every $n\in\mathbb N_0$, where $R_{n,s}$ is the set of all sequences $(r_0,...,r_n)$ of nonnegative integers such that $r_0+...+r_n=s$, $r_1+...+nr_n=n$. This proves (a).
\item[(b)] $r(g)>0$:\\
Let us extract from the sequence $(f_k)$ two distinct infinite subsequences, if possible -- let $(f_k^{I})_{k\in\mathbb Z_+}$ be a subsequence of $(f_k)$ consisting of all such $f_k$ that $a_{0,k}\neq 0$ and let $(f_k^{II})_{k\in\mathbb Z_+}$ be a subsequence of $(f_k)$ consisting of all $f_k$ such that $a_{0,k}=0$. If this is not possible (one of the subsequences $f_k^{I,II}$ would be only finite or would not exist at all), it is sufficient to consider the one existing infinite subsequence and repeat the below reasoning. Then one can prove analogously to the proof of a) that $\lim\limits_{k\rightarrow\infty}g\circ f_k^{II}=g\circ f$. We will now prove that $\lim\limits_{k\rightarrow\infty}g\circ f_k^{I}=g\circ f$. Denote $f_k^{I}=\sum\limits_{n=0}^{\infty}a_{n,k}^{I}z^n$ for every $k\in\mathbb Z_+$. We have $g\circ f_k^{I}=\sum\limits_{n=0}^{\infty}c_{n,k}^{I}z^n$, where
\[
c_{n,k}^{I}=\sum\limits_{s=0}^{\infty}b_sa_{n,k}^{I\,(s)}=\sum\limits_{s=0}^{\infty}\sum\limits_{R_{n,s}}\frac{b_ss!(a_{0,k}^I)^{r_0}\ldots (a_{n,k}^I)^{r_n}}{r_0!\ldots r_n!}=\sum\limits_{s=0}^{\infty}\sum\limits_{r_1+...+nr_n=n}\frac{b_ss!(a_{0,k}^I)^{s-r_1-...-r_n}\ldots (a_{n,k}^I)^{r_n}}{(s-r_1-...-r_n)!\ldots r_n!},
\]
where we admit $\frac{s!}{(s-r_1-...-r_n)!}=0$ for $s<r_1+...+r_n$. \\
First assume $\mbox{deg}(f)>0$. Then we can assume the degree of all $f_k^I$ is greater than $0$ for sufficiently large $k$. It is easy to check that for every $t\in\mathbb N_0$, $g_t(z):=\sum\limits_{s=t}^{\infty}\frac{s!}{(s-t)!}b_sz^s$ is a power series with radius of convergence $r(g)$. For all $z_0\in\mathbb C$, the convergence of $f_t(z_0)$ for every $t\in\mathbb N_0$ is equivalent to the convergence of $g^{(t)}(z_0)$ for every $t\in\mathbb N_0$. Therefore, by the General Composition Theorem (Thm. \ref{gct}), the series $\sum\limits_{s=0}^{\infty}\frac{b_ss!}{(s-r_1-...-r_n)!}(a_{0,k}^I)^s$ is convergent for every $r_1,...,r_n\in\mathbb N_0$ (since $g\circ f_k^{I}$ exists for every $k$), and
\[
c_{n,k}^{I}=\sum\limits_{r_1+...+nr_n=n}\left(\frac{(a_{1,k}^I)^{r_1}\ldots (a_{n,k}^I)^{r_n}}{(a_{0,k})^{r_1+...+r_n}r_1!\ldots r_n!}\sum\limits_{s=0}^{\infty}\frac{b_ss!}{(s-r_1-...-r_n)!}(a_{0,k}^I)^s\right).
\]
Now
\begin{itemize}
\item if $|a_0|<r(g)$, then $\sum\limits_{s=0}^{\infty}\frac{b_ss!}{(s-r_1-...-r_n)!}(a_{0,k}^I)^s\stackrel{k\rightarrow\infty}{\rightarrow}\sum\limits_{s=0}^{\infty}\frac{b_ss!}{(s-r_1-...-r_n)!}(a_{0})^s$ (the sum of a power series is a continuous function in every interior point of its ball of convergence), so $c_{n,k}^{I}\stackrel{k\rightarrow\infty}{\rightarrow}c_n$ for every $n\in\mathbb N_0$,
\item if $r(g)<+\infty$ and $|a_0|=r(g)$, then by Thm. \ref{l9.2.5}, $\sum\limits_{s=0}^{\infty}\frac{b_ss!(a_{0,k}^I)^s}{(s-r_1-...-r_n)!}\stackrel{k\rightarrow\infty}{\rightarrow}\sum\limits_{s=0}^{\infty}\frac{b_ss!a_{0}^s}{(s-r_1-...-r_n)!}$, so $c_{n,k}^{I}\stackrel{k\rightarrow\infty}{\rightarrow}c_n$ for every $n\in\mathbb N_0$.
\end{itemize}
Now assume $\mbox{deg}(f)=0$ and $|a_0|<r(g)$. Then there exists such $\varepsilon>0$ and $K\in\mathbb N$ that for every $k>K$, $|a_{0,k}^I|<r(g)-\varepsilon$. Therefore for sufficiently large $k$, it can be proven analogously to the case $\mbox{deg}(f)>0$ that $\sum\limits_{s=0}^{\infty}\frac{b_ss!}{(s-r_1-...-r_n)!}(a_{0,k}^I)^s$ is convergent for every $r_1,...,r_n\in\mathbb N_0$ and
\begin{eqnarray*}
\lim\limits_{k\rightarrow\infty}c_{n,k}^I&=&\lim\limits_{k\rightarrow\infty}\sum\limits_{r_1+...+nr_n=n}\frac{(a_{1,k}^I)^{r_1}\ldots (a_{n,k}^I)^{r_n}}{(a_{0,k})^{r_1+...+r_n}r_1!\ldots r_n!}\sum\limits_{s=0}^{\infty}\frac{b_ss!}{(s-r_1-...-r_n)!}(a_{0,k}^I)^s\\&=&\sum\limits_{r_1+...+nr_n=n}\frac{a_{1}^{r_1}\ldots a_n^{r_n}}{(a_{0})^{r_1+...+r_n}r_1!\ldots r_n!}\sum\limits_{s=0}^{\infty}\frac{b_ss!}{(s-r_1-...-r_n)!}a_0^s=c_n
\end{eqnarray*}
for every $n\in\mathbb N_0$, which completes the proof.
\end{enumerate}
\end{proof}


\begin{remark}
A more general mapping $D\ni (g,f)\mapsto g\circ f\in \mathbb X(\mathbb C)$, where $D$ is a subset of $\mathbb X(\mathbb C)\times \mathbb X(\mathbb C)$ such that all the compositions exist is not continuous in general -- indeed, let $f_k=2^k\sum\limits_{n=k}^{\infty}z^n$ for every $k\in\mathbb N$ and let $g=b_0+b_1z+...$ be a formal power series with $b_0=\frac12$. Then $|b_0|<r(f_k)$, so $f_k\circ g$ exists for every $k\in\mathbb N$ and $f_k\stackrel{k\rightarrow\infty}{\rightarrow}0$. However, it is easy to check that for every $k\in\mathbb N$, the $0$th coefficient of $f_k\circ g$ is equal to $f_k(\frac12)=2$.
\end{remark}
To establish a necessary and sufficient condition for the existence of left composition inverse of a formal power series we will need two lemmas, which -- however well-known -- we provide with proofs for the convenience of the reader.
\begin{lemma}\label{rek}
Let $f(z)=a_0+a_1z+a_2z^2+\ldots$, $a_0=0$, $a_1\neq 0$ be a formal power series and denote $f^k(z)=a_k^{(k)}z^k+a_{k+1}^{(k)}z^{k+1}+...$ for every $k\in\mathbb N_0$. Then $g(z)=b_0+b_1z+b_2z^2+...$, where
\[
\left\{\begin{array}{l}
b_0=0,\\
b_1=\frac{1}{a_1},\\
b_n=-\frac{b_1a_n^{(1)}+...+b_{n-1}a_n^{(n-1)}}{a_1^n}\\
\end{array}\right.
\]
satisfies $g\circ f(z)=z$. Moreover, $f$ possesses no other composition inverses.
\end{lemma}

\begin{proof}
(Cf. \cite{g21}, proof of Thm. 1.5.9.) Denote $g\circ f(z)=c_0+c_1z+c_2z^2+...$. We have $c_0=a_0b_0^{(0)}=0$, $c_1=b_0a_1^{(0)}+b_1a_1^{(1)}=1$ and
\[
c_n=\sum\limits_{k=0}^nb_ka_n^{(k)}=b_na_n^{(n)}+\sum\limits_{k=1}^{n-1}b_ka_n^{(k)}=b_na_1^n+\sum\limits_{k=1}^{n-1}b_ka_n^{(k)}=0,\,\,n>1.
\]
Now, let $h(z)=\sum\limits_{n=0}^{\infty}h_nz^n\in\mathbb X(\mathbb C)$ be a formal power series such that $f\in\mathbb X_h$ and $h\circ f(z)=z$. Then $h_0a_0^{(0)}=h_0=0$, $h_0a_1^{(0)}+h_1a_1^{(1)}=h_1a_1=1$ and for every $n>1$, $\sum\limits_{k=0}^nh_ka_n^{(k)}=h_na_1^n+\sum\limits_{k=1}^{n-1}h_ka_n^{(k)}=0$, which proves $h=g$.
\end{proof}

\begin{lemma}\label{comp}
Let $f(z)=\sum\limits_{n=0}^{\infty}a_nz^n$, $g(z)=\sum\limits_{n=0}^{\infty}b_nz^n\in\mathbb X(\mathbb C)$. If $f\in\mathbb X_g$, then $g\circ f=g_D\circ(f-a_0)$, where
$g_D(z)=\sum\limits_{n=0}^{\infty}\frac{g^{(n)}(a_0)}{n!}z^n$.
\end{lemma}

\begin{proof}
(Cf. \cite{g21}, proof of Thm. 5.4.1.) Since $g\circ f$ exists, $g^{(k)}(a_0)=\sum\limits_{n=k}^{\infty}\frac{n!}{(n-k)!}b_na_0^{n-k}$ exists for every $k\in\mathbb N_0$. Denote $f-a_0(z)=A_0+A_1z+...$ . We have
\[
[z^0]g_D\circ (f-a_0)=g(a_0)A_0^{(0)}=\sum\limits_{n=0}^{\infty}b_na_0^n=\sum\limits_{n=0}^{\infty}b_na_0^{(n)}=[z^0]g\circ f
\]
and, for every $k>0$,
\begin{eqnarray*}
[z^k]g_D\circ (f-a_0)&=&\sum\limits_{j=1}^k\frac{g^{(j)}(a_0)}{j!}A_k^{(j)}=\sum\limits_{j=1}^kA_k^{(j)}\left(\sum\limits_{n=j}^{k-1}\binom{n}{j}b_na_0^{n-j}+\sum\limits_{n=k}^{\infty}\binom{n}{j}b_na_0^{n-j}\right)\\&=&\sum\limits_{j=1}^kA_k^{(j)}\left(\sum\limits_{n=j}^{k-1}\binom{n}{j}b_na_0^{n-j}\right)+\sum\limits_{n=k}^{\infty}b_n\left(\sum\limits_{j=1}^k\binom{n}{j}a_0^{n-j}A_k^{(j)}\right)\\&=&\sum\limits_{n=1}^{k-1}b_n\sum\limits_{j=1}^n\binom{n}{j}A_k^{(j)}a_0^{n-j}+\sum\limits_{n=k}^{\infty}b_na_k^{(n)}=\sum\limits_{n=1}^{\infty}b_na_k^{(n)}=\sum\limits_{n=0}^{\infty}b_na_k^{(n)}=[z^k]g\circ f,
\end{eqnarray*}
which completes the proof.
\end{proof}

\begin{definition}
Let $A\in\mathbb C$. We define the generalized Pascal matrix $P(A)$ by the formula
\[
P(A)=\left[\begin{array}{lllllll}
1 & A & A^2 & \ldots & A^n & \ldots \\
0 & 1 & 2A & \ldots & nA^{n-1} & \ldots \\
0 & 0 & 1 & \ldots & \binom{n}{2}A^{n-2} & \ldots \\
\vdots & \vdots & \vdots & \ddots & \vdots & \ddots \\
0 & 0 & \ldots & \ldots & \binom{n}{k}A^{n-k} & \ldots \\
\vdots & \vdots & \ddots & \ddots & \vdots & \ddots \\
\end{array}\right].
\]
\end{definition}
\begin{theorem}\label{inv}
Let $f(z)=\sum\limits_{n=0}^{\infty}a_nz^n\in\mathbb X(\mathbb C)$. Then there exists such $g\in\mathbb X(\mathbb C)$ that $g\circ f(z)=z$, if and only if $a_1\neq 0$ and there exists such $\tilde{g}\in\mathbb X(\mathbb C)$ that $\tilde{g}\circ(z+a_0)=(f-a_0)^{[-1]}$. If such a $\tilde{g}$ exists, then $g=\tilde{g}$.\\
Equivalently: there exists such $g\in\mathbb X(\mathbb C)$ that $g\circ f(z)=z$, if and only if $a_1\neq 0$ and the infinite system $P(a_0)[b_0\,b_1\,...]^T=[c_0\,c_1\,...]^T$ ($c_n:=[z^n](f-a_0)^{[-1]}$) possesses a solution. If such a $g$ exists, then $g=b_0+b_1z+...$.
\end{theorem}

\begin{proof}
Let first $a_1=0$. Then for every $k\in\mathbb N_0$, $a_1^{(k)}=0$, so for any $g$ such that $f\in\mathbb X_g$, $[z^1]g\circ f=0$, so $g\circ f\neq z$.\\
Assume there exists such $g(z)=b_0+b_1z+b_2z^2+...$ that $g\circ f(z)=z$. Then $a_1\neq 0$ and, by Lemma \ref{comp}, $g_D\circ(f-a_0)(z)=z$, so $(f-a_0)^{[-1]}(z)=g_D=\sum\limits_{k=0}^{\infty}\left(\sum\limits_{n=k}^{\infty}\binom{n}{k}b_na_0^{n-k}\right)z^k$ (since $(f-a_0)^{[-1]}$ exists and is unique -- see Lemma \ref{rek}). Denote $(f-a_0)^{[-1]}(z)=c_0+c_1z+c_2z^2+...$. We have
\[
\left[\begin{array}{lllllll}
1 & a_0 & a_0^2 & \ldots & a_0^n & \ldots \\
0 & 1 & 2a_0& \ldots & na_0^{n-1} & \ldots \\
0 & 0 & 1 & \ldots & \binom{n}{2}a_0^{n-2} & \ldots \\
\vdots & \vdots & \vdots & \ddots & \vdots & \ddots \\
0 & 0 & \ldots & \ldots & \binom{n}{k}a_0^{n-k} & \ldots \\
\vdots & \vdots & \ddots & \ddots & \vdots & \ddots \\
\end{array}\right]\left[\begin{array}{l}
b_0\\
b_1\\
b_2\\
\vdots\\
b_k\\
\vdots\\
\end{array}\right]=P(a_0)\left[\begin{array}{l}
b_0\\
b_1\\
b_2\\
\vdots\\
b_k\\
\vdots\\
\end{array}\right]=\left[\begin{array}{l}
c_0\\
c_1\\
c_2\\
\vdots\\
c_k\\
\vdots\\
\end{array}\right].
\]
It is also easy to check that $P(a_0)=C^T_{z+a_0}$ (the composition matrix of the formal series $z+a_0$, see Preliminaries), which completes this part of proof.\\
Now assume $a_1\neq 0$ and $\tilde{g}$ exists, 
or, in other words, the infinite system $P(a_0)[b_0\,b_1\,...\,b_k\,...]^T=[c_0\,c_1\,...\,c_k\,...]^T$ (with unknowns $b_i$) possesses a solution. It is easy to check using the Multinomial Theorem that for every formal power series $f=a_0+a_1z+...$,
\begin{eqnarray*}
[C_f^T\,|\,[0\,1\,0...]^T\,]
&=&\left[\begin{array}{llllll|l}
a_0^{(0)} & a_0^{(1)} & a_0^{(2)} & \ldots & a_0^{(n)} & \ldots & 0 \\
a_1^{(0)} & a_1^{(1)} & a_1^{(2)} & \ldots & a_1^{(n)} & \ldots & 1\\
a_2^{(0)} & a_2^{(1)} & a_2^{(2)} & \ldots & a_2^{(n)} & \ldots & 0\\
\vdots & \vdots & \vdots & \ddots & \vdots & \ddots &\vdots\\
a_k^{(0)} & a_k^{(1)} & \ldots & \ldots & a_k^{(n)} & \ldots &0\\
\vdots & \vdots & \ddots & \ddots & \vdots & \ddots &\vdots\\
\end{array}\right]\\
&\cong&\left[\begin{array}{llllll|l}
1 & a_0 & a_0^2 & \ldots & a_0^n & \ldots & q_0 \\
0 & 1 & 2a_0 & \ldots & na_0^{n-1} & \ldots & q_1\\
0 & 0 & 1 & \ldots & \binom{n}{2}a_0^{n-2} & \ldots & q_2\\
\vdots & \vdots & \vdots & \ddots & \vdots & \ddots &\vdots\\
0 & 0 & \ldots & \ldots & \binom{n}{k}a_0^{n-k} & \ldots &q_k\\
\vdots & \vdots & \ddots & \ddots & \vdots & \ddots &\vdots\\
\end{array}\right]
=[P(a_0)\,|\,[q_0...q_k...]^T\,],
\end{eqnarray*}
where
\begin{enumerate}
\item "$\cong$" means "equivalent with respect to row reduction (Gaussian elimination)" and $C_f$ denotes the composition matrix of $f$,
\item one needs a finite number of operations for each row to get properly reduced; also, this property holds for the ''inverted'' operation, ''$P(a_0)\rightarrow C_f^T$'' as well;
\item $q_0,q_1,q_2,...$ are some complex coefficients that do not depend on $a_0$.
\end{enumerate}
We will now find the coefficients $q_n$. Their values are independent from $a_0$, so one can assume that $a_0=0$. Then $P(a_0)=I$ (the unit matrix), so $[q_0...q_k...]^T$ is the solution to the equation $C^T_{f-a_0}[b_0...b_k...]^T=[0\,1\,0...0...]^T$ (with unknowns $b_i$). Therefore, by Prop. \ref{matrixcomp}, $(q_0+q_1z+...)\circ (f-a_0) =z$, so $\sum\limits_{n=0}^{\infty}q_nz^n=(f-a_0)^{[-1]}$ (for nonunit formal power series, the left and right composition inverse are always equal if they exist, see e.g. \cite{g21}, Lemma 1.5.5) so $(q_n)_{n\in\mathbb N_0}=(c_n)_{n\in\mathbb N_0}$ (Lemma \ref{rek}). It follows that if the system $P(a_0)[b_0...b_k...]^T=[c_0\,c_1\,...\,c_k\,...]^T$ possesses a solution, then $C_f^T[b_0...b_k...]^T=[0\,1\,0\,...]^T$ possesses a solution, which, by Prop. \ref{matrixcomp}, completes the proof.
\end{proof}
Let us emphasize that infinite linear systems like $P(a_0)[b_0...b_k...]^T=[c_0\,c_1\,...\,c_k\,...]^T$, which appears in the above theorem, have been extensively studied for decades and there are some well-known results of finding their solutions (and when they exist) -- see e.g. \cite{fed}. Also, we have the following
\begin{corollary}
For every $f\in\mathbb X(\mathbb C)$, its left composition inverse is unique if it exists.
\end{corollary}
\begin{proof}
Let $g=b_0+b_1z+...$, $g\circ f=z$ and denote $(f-a_0)^{[-1]}=c_0+c_1z+...$ . Then $P(a_0)[b_0\,b_1\,...]^T=[c_0\,c_1\,...]^T$. Assume there exists such a sequence $(d_n)_{n\in\mathbb N_0}$ that $P(a_0)[d_0\,d_1\,...]^T=[0\,0\,...]^T$; then the sequence $(e_n)_{n\in\mathbb N_0}$, $e_n=(n+1)d_{n+1}$ also satisfies this homogenous infinite system. Therefore by e.g. \cite{fed}, Thm. 5.2 there exists such $c\in\mathbb C$ that $e_i=cd_i$ for every $i\in\mathbb N_0$, so $d_{i}=\frac{cd_0}{i!}$ -- but this sequence is not a solution to the above infinite system, which completes the proof.
\end{proof}

Let us now denote as $\mathbb C^{\mathbb N}$ the set of all complex sequences $(a_n)_{n\in\mathbb N_0}$, with the family of seminorms $\|(a_n)\|_k=|a_k|$ ($k\in\mathbb N_0$), so that it forms a Fr\'{e}chet space. It is obvious that $(\mathbb X(\mathbb C),A)$, where $\phi:\mathbb X(\mathbb C)\ni \sum\limits_{n=0}^{\infty}a_nz^n\to (a_0,a_1,a_2,...)\in \mathbb C^{\mathbb N}$ and $A$ is the maximal atlas compatible with $\left\{(\mathbb X(\mathbb C),\phi)\right\}$, is a Fr\'{e}chet manifold modeled on $\mathbb C^{\mathbb N}$. This setting allows us to formulate the following
\begin{theorem}\label{smooth}
Let $g\in\mathbb X(\mathbb C)$ and let $T_g:\mathbb X_g\ni f\mapsto g\circ f\in\mathbb X(\mathbb C)$. Then
\begin{enumerate}
\item[a)] if $r(g)=0$, then $T_g$ is smooth;
\item[b)] if $r(g)>0$, then $T_g$ is smooth in every point $f=\sum\limits_{n=0}^{\infty}a_nz^n$ such that (1) $|a_0|<r(g)$ or (2) $r(g)<+\infty$, $|a_0|=r(g)$ and $\mbox{deg}(f)>0$.
\end{enumerate}
\end{theorem}
Before presenting the proof of the above result, a clarification regarding the definition of smoothness is necessary since the domain of $T_g$ ($\mathbb X_g$) is not a priori an open subset of $\mathbb X(\mathbb C)$. Similarly to Theorem \ref{thm1}, we endow $\mathbb X_g$ with a subspace topology and observe that $\mathbb X_g=\left\{a_0+a_1z+...:a_0\in D\right\}$, where $D$ is either $\left\{0\right\}$, an open disk centered at 0, a closed disk centered at 0 or the whole complex plane. We then define after e.g. \cite{monastir} the derivative of $T_g$ at point $w$ in direction $k$ (if it exists) by the formula $\lim\limits_{t\rightarrow 0}\frac{T_g(w+tk)-T_g(w)}{t}$, where, if $w$ lies on the boundary of $\mathbb X_g$, we restrict this limit $t\rightarrow 0$ to limit over sequences $(t_n)$ converging to $0$ such that for every $n$, $w+t_nk\in\mathbb X_g$ instead of all possible sequences $(t_n)$ converging to $0$. The notion of smooth ($C^{\infty}$) maps is then analogous to \cite{monastir}.\\

Let us now present the proof of Theorem \ref{smooth}: 
\begin{proof}
Denote $g=b_0+b_1z+...$. Claim a) is obvious -- if $r(g)=0$, than all series $f\in\mathbb X_g$ are nonunits, so for every $n\in\mathbb N_0$, $[z^n]T_g(f)$ is a polynomial function of $b_0,...,b_n$ and $[z^1]f,...,[z^n]f$. Now let $r(g)>0$, $w=w_0+w_1z+...\in\mathbb X_g$, and $k=k_0+k_1z+...\in\mathbb X(\mathbb C)$. Let us first prove the following simple fact, which we will call Step I: 
\\
{\bf Step I.}
{\it Let $A,B\in\mathbb N_0$. Then for $t$ sufficiently close to 0 (in the sense clarified above the proof) the series
\[
F_{A,B}(t):=\sum\limits_{s=A+\max(2,B)}^{\infty}b_s\sum\limits_{l=A}^{s-\max(2,B)}t^{s-2-l}\frac{s!}{(l-A)!(s-l-B)!}w_0^{l-A}k_0^{s-l-B}
\] is convergent and $\lim\limits_{t\rightarrow 0}F_{A,B}(t)=0$.}\\

\textit{Proof of Step I.} Let us divide the proof into three cases:
\begin{enumerate}
\item $B\geq 2$:\\
We have
\begin{eqnarray*}
&&
\sum\limits_{s=A+B}^{\infty}b_s\sum\limits_{l=A}^{s-B}t^{s-1-l}\frac{s!}{(l-A)!(s-l-B)!}w_0^{l-A}k_0^{s-l-B}\\
&&
=t^{B-1}\sum\limits_{s=A+B}^{\infty}\frac{s!}{(s-A-B)!}b_s\sum\limits_{l=A}^{s-B}\binom{s-A-B}{l-A}w_0^{l-A}k_0^{s-l-B}t^{s-l-B}\\
&&
=t^{B-1}\sum\limits_{s=0}^{\infty}\frac{(s+A+B)!}{s!}b_{s+A+B}(w_0+k_0t)^s.
\end{eqnarray*}
The last series is convergent since $g\circ (w+tk)$ exists (Thm. \ref{gct}), and therefore $F_{A,B}(t)$ exists. Now, obviously $t^{B-1}\rightarrow 0$ and
\[
\lim\limits_{t\rightarrow 0}\sum\limits_{s=0}^{\infty}\frac{(s+A+B)!}{s!}b_{s+A+B}(w_0+k_0t)^s=\sum\limits_{s=0}^{\infty}\frac{(s+A+B)!}{s!}b_{s+A+B}w_0^s<+\infty,
\]
since a) if $|w_0|<r(g)$, then the last equality is a consequence of continuity of functions given by power series on their (open) disk of convergence b) if $r(g)<+\infty$, $|a_0|=r(g)$ and $\mbox{deg}(f)>0$, then the last equality results from Thm. \ref{l9.2.5}.
\item $B=1$:\\
Now we have
\begin{eqnarray*}
&&
\sum\limits_{s=A+1}^{\infty}b_s\sum\limits_{l=A}^{s-2}t^{s-1-l}\frac{s!}{(l-A)!(s-l-1)!}w_0^{l-A}k_0^{s-l-1}\\
&&
=\sum\limits_{s=A+1}^{\infty}\frac{s!}{(s-A-1)!}b_s\left((w_0+k_0t)^{s-A-1}-w_0^{s-A-1}\right)\\
&&
=\sum\limits_{s=0}^{\infty}\frac{(s+A+1)!}{s!}b_{s+A+1}\left((w_0+k_0t)^s-w_0^s\right)
\end{eqnarray*}
The limit of the above expression as $t\rightarrow 0$ is equal to $0$ by an argument analogous as in (1).
\item $B=0$:\\
In this case
\begin{align*}
&\sum\limits_{s=A}^{\infty}\frac{s!}{(s-A)!}b_s\sum\limits_{l=A}^{s-2}t^{s-1-l}\binom{s-A}{l-A}w_0^{l-A}k_0^{s-l}\\
&=\sum\limits_{s=0}^{\infty}\frac{(s+A)!}{s!}b_{s+A}\left(\frac{(w_0+tk_0)^s-w_0^s}{t}-sk_0w_0^{s-1}\right)\\
&=k_0\left(\frac{g^{(A)}(w_0+k_0t)-g^{(A)}(w_0)}{k_0t}-g^{(A+1)}(w_0)\right)
\end{align*}
which also tends to $0$ when $t\rightarrow 0$ (Thm. \ref{l9.2.5} and basic properties of differentiation of power series), which completes the proof of Step I.
\end{enumerate}
Now, returning to the main part of the proof, we have
\begin{align*}
\lim\limits_{t\rightarrow 0}\frac{g\circ (w+tk)-g\circ w}{t}&=\sum\limits_{n=0}^{\infty}\lim\limits_{t\rightarrow 0}\left(\sum\limits_{s\in\mathbb N_0}b_s\frac{[z^n](w+tk)^s-[z^n]w^s}{t}\right)z^n\\&=\sum\limits_{n=0}^{\infty}\lim\limits_{t\rightarrow 0}\left(\sum\limits_{s\in\mathbb N}b_s\sum\limits_{l=0}^{s-1}\binom{s}{l}t^{s-l-1}[z^n](w^lk^{s-l})\right)z^n\\&=\sum\limits_{n=0}^{\infty}\lim\limits_{t\rightarrow 0}\left(\sum\limits_{s\in\mathbb N}b_s\left(s[z^n](w^{s-1}k)+\sum\limits_{l=0}^{s-2}\binom{s}{l}t^{s-l-1}[z^n](w^lk^{s-l})\right)\right)z^n,
\end{align*}
where -- here and later in the proof -- we denote $\sum\limits_{l=0}^{N}(\ldots):=0$ for $N<0$. Now, see that $g'\circ w\in\mathbb X(\mathbb C)$ (\cite{g21}, Lemma 5.5.2. if deg$(w)\neq 0$; obvious if deg$(w)=0$ and $|w_0|<r(g)$) and
\begin{multline*}
(g'\circ w)k=
\sum\limits_{n=0}^{\infty}\left(\sum\limits_{s=0}^{\infty}(s+1)b_{s+1}[z^n]w^s\right)z^n\, k=
\sum\limits_{s=0}^{\infty}\left((s+1)b_{s+1}\sum\limits_{n=0}^{\infty}[z^n]w^s\right)z^n\, k\\=
\sum\limits_{s=0}^{\infty}\left((s+1)b_{s+1}w^s\right)k=
\sum\limits_{s=0}^{\infty}\left((s+1)b_{s+1}\sum\limits_{n=0}^{\infty}[z^n](w^sk)z^n\right)=
\sum\limits_{n=0}^{\infty}\left(\sum\limits_{s=1}^{\infty}sb_{s}[z^n](w^{s-1}k)\right)z^n
\end{multline*}
so $\sum\limits_{s=1}^{\infty}sb_{s}(w^{s-1}k)_n$ exists for every $n\in\mathbb N_0$. Notice that the above substituting $\sum\limits_{s=0}^{\infty}\sum\limits_{n=0}^{\infty}$ for $\sum\limits_{n=0}^{\infty}\sum\limits_{s=0}^{\infty}$ and vice versa is not changing the order of summation in a double series -- the summation over $n$ is just a way to denote a formal power series, and what we did was using the topology defined earlier in this article on $\mathbb X(\mathbb C)$ to sum an infinite sequence of formal power series ''term by term''. We also used a simple fact that the multiplication of formal power series is continuous with respect to that topology (for every $n\in\mathbb N$, the mappings $\mathbb C\times\ldots\times\mathbb C\ni(a_0,\ldots,a_n,b_0,\ldots,b_n)\mapsto \sum\limits_{k=0}^na_kb_{n-k}$ are obviously continuous).\\
Therefore $\sum\limits_{s=2}^{\infty}b_s\sum\limits_{l=0}^{s-2}\binom{s}{l}t^{s-l-1}[z^n](w^lk^{s-l})\in\mathbb C$ for every $n\in\mathbb N$ and $t$ sufficiently close to $0$ (in the sense clarified above this proof) and
\[
\lim\limits_{t\rightarrow 0}\frac{g\circ (w+tk)-g\circ w}{t}=(g'\circ w)k+\lim\limits_{t\rightarrow 0}\sum\limits_{n=0}^{\infty}\left(\sum\limits_{s=2}^{\infty}b_s\sum\limits_{l=0}^{s-2}\binom{s}{l}t^{s-l-1}[z^n](w^lk^{s-l})\right)z^n.
\]

Now, using the Multinomial Theorem, we have, for every $n\in\mathbb N_0$,
\begin{align*}
&\sum\limits_{s=2}^{\infty}b_s\sum\limits_{l=0}^{s-2}\binom{s}{l}t^{s-l-1}[z^n](w^lk^{s-l})\\
&=\sum\limits_{s=2}^{\infty}\sum\limits_{l=0}^{s-2}b_st^{s-1-l}\binom{s}{l}\sum\limits_{T=0}^n
\left[\left(\sum\limits_{r(T)}\frac{l!}{r_1!...r_T!(l-r_1-...-r_T)!}w_1^{r_1}...w_T^{r_T}w_0^{l-r_1-...-r_T}\right)\right.\\
&\hspace*{1cm}\left.
\left(\sum\limits_{R(n-T)}\frac{(s-l)!}{R_1!...R_{n-T}!(s-l-R_1-...-R_{n-T})!}k_1^{R_1}...k_{n-T}^{R_{n-T}}k_0^{s-l-R_1-...-R_{n-T}}\right)\right]\\
&=\sum\limits_{T=0}^nb_s\sum\limits_{r(T),R(n-T)}\frac{w_1^{r_1}...w_T^{r_T}k_1^{R_1}...k_{n-T}^{R_{n-T}}}{r_1!...r_T!R_1!...R_{n-T}!}F_{r_1+...+r_T,R_1+...+R_{n-T}}(t)\stackrel{t\rightarrow 0}{\longrightarrow}0
\end{align*}
(see Step I), where $\frac{l!}{(l-q)!}:=0$ if $q>l$ and $r(T)=\left\{(r_1,...,r_T)\in\mathbb N_0^T:r_1+...+Tr_T=T\right\}$ (and we define $R_{n-T}$ analogously). This proves that the derivative of $T_g$ in considered points in any direction $k$ exists and is equal to $(g'\circ w)k$, which is continuous with respect to $k$ and $w$. The proof concerning higher order derivatives of $T_g$ is analogous.
\end{proof}
\begin{remark}
The claims of Theorems \ref{thm1}, \ref{smooth} also hold if deg$(f)=0$, $|[z^0]f|=r(g)$, provided that for some $a\in\mathbb C$, $|a|=r(g)$, $g^{(k)}(a)\in\mathbb C$ for every $k\in\mathbb N_0$ (Thm. \ref{l9.2.5}).
\end{remark}
By calculations analogous to the ones in the above proof, one can conclude the following
\begin{corollary}
("Taylor formula" for the general composition of formal power series) Let $g,f,k\in\mathbb X(\mathbb C)$, where if $r(g)\neq 0$, then we assume $|[z^0]f|<r(g)$ or $r(g)<\infty$, $|[z^0]f|=r(g)$ and $\mbox{deg}(f)>0$. Then
\[
\lim\limits_{t\rightarrow 0}\frac{g\circ (f+tk)-\left(g\circ f+t(g'\circ f)k+...+\frac{t^n}{n!}(g^{(n)}\circ f)k^n\right)}{t^n}=0,
\]
where $\lim\limits_{t\rightarrow 0}$ has the same meaning as in Theorem \ref{smooth}.
\end{corollary}

The findings presented in this section bring us closer to the possibility of employing methods of differential geometry and Lie theory to consider the general composition of formal power series. This issue will be a topic of our further exploration in the future. In this article, we are now going to analyze some properties of the Fr\'{e}chet-Lie group structures on the set of nonunit formal power series.

\section{Fr\'{e}chet-Lie group structures on some families of formal power series}
The so called substitution group $\xi(R)=\left\{z+a_2z^2+a_3z^3+...:a_2,a_3,...\in R\right\}$ of formal power series over a commutative ring $R$ has been first introduced by Jennings in 1954 \cite{jennings}. It is a group of all nonunit formal power series $f$ with $a_1=1$ and the composition of formal power series $\circ$ as the group action. We will denote this group by $\mathbb X^0_1(R)$ in this paper; we will also assume $R=\mathbb R$ or $R=\mathbb C$ with its natural topology; then the topology on $\mathbb X^0_1(R)$ is equivalent to the product topology on $\mathbb C^{\infty}$ ($\mathbb R^{\infty}$). Let $\phi:\mathbb X_1^0(\mathbb C)\ni z+\sum\limits_{n=2}^{\infty}a_nz^n\to (a_2,a_3,...)\in \mathbb C^{\mathbb N}\,(\mathbb R^{\mathbb N})$ (the Fr\'{e}chet space of all complex (real) sequences) and let $A$ be the maximal atlas compatible with $\left\{(\mathbb X_1^0(\mathbb C),\phi)\right\}$. Then it is well-known that $((\mathbb X_1^0(\mathbb C),\circ,\tau),A)$ is a Fr\'{e}chet-Lie group modeled on $\mathbb C^{\mathbb N}$ ($\mathbb R^{\mathbb N})$. The Lie algebra of this Fr\'{e}chet-Lie group is the algebra of all formal series of the form $a_2z^2+a_3z^3+...$ with Lie bracket $[f,g]=fg'-f'g$.\\
A broader group $\mathcal{A}(\mathbb C)$ -- a group of all nonunit formal power series $f$ with $a_1\neq 0$ and $\circ$ as the group action (here we will denote it as $\mathbb X_z^0(\mathbb C)$) -- has also been mentioned in literature, mainly in context of finding conjugacy classes of formal power series and the Schröder's equation \cite{babenko}. Here we give a systematic description of the group $\mathbb X_z^0(\mathbb C)$ as a Fr\'{e}chet-Lie group and some of its properties. 
\begin{remark}
See that an other approach than the one presented below would be to immediately associate $\mathbb X_z^0(\mathbb C)$ with the semidirect product $\mathbb X_1^0(\mathbb C)\rtimes_\phi \mathbb C^*$, where the homomorphism $\phi$ is defined as $\phi:\mathbb C^*\ni a\mapsto (g\mapsto ag(a^{-1}z))\in\textrm{Aut}(\mathbb X_1^0(\mathbb C))$. Although it would be a little briefer way to obtain the below results, we provide here, for the convenience of the reader, as well as for the sake of future investigations of more complex structures like the Riordan groups, an explicit definition and systematic derivation of the Lie group and corresponding Lie algebra properties without employing the semidirect product structure from the beginning.
\end{remark}
First, let us introduce the following simple
\begin{proposition}\label{proplie1}
Let $\phi:\mathbb X_z^0(\mathbb C)\ni \sum\limits_{n=1}^{\infty}a_nz^n\to (a_1,a_2,...)\in \mathbb C^{\mathbb N}$ and let $A$ be the maximal atlas compatible with $\left\{(\mathbb X_z^0(\mathbb C),\phi)\right\}$. Then $((\mathbb X_z^0(\mathbb C),\circ,\tau),A)$ is a Fr\'{e}chet-Lie group modeled on $\mathbb C^{\mathbb N}$.
\end{proposition}
\begin{proof}
It is obvious that $s_z:=\phi(\mathbb X_z^0(\mathbb C))=\{(a_n)_{n\in\mathbb N}\in \mathbb C^{\mathbb N}:a_1\neq 0\}$ is an open subset of $\mathbb C^{\mathbb N}$ and that $\phi$ is a homeomorphism from $\mathbb X_z^0(\mathbb C)$ to $s_z$. Moreover, the mappings $(g,f)\mapsto g\circ f$, $g\mapsto g^{[-1]}$ are smooth -- for example, let $T:g\to g^{[-1]}$ and let $a=(a_1,a_2,...)$, $a_1\neq 0$. By Lemma \ref{rek},
\[
(\phi\circ T\circ\phi^{-1})(a)=(W_1(a_1),W_2(a_1,a_2),...,W_n(a_1,...,a_n),...),
\]
where $W_n$ are some polynomial functions of $a_2,...,a_n$, and rational functions of $a_1$. Therefore for every $n\in\mathbb N$, $W_n$ is holomorphic on $(\mathbb C\setminus\left\{0\right\})\times \mathbb C ^{n-1}$.
\end{proof}
\begin{proposition}
The Lie algebra corresponding to the Fr\'{e}chet-Lie group $\mathbb X_z^0(\mathbb C)$ can be identified with $(\mathbb X^0(\mathbb C),[\cdot])$ (the set of all nonunit formal power series), where the Lie bracket $[\cdot]$ is given by a formula the same as for $\mathbb X_1^0(\mathbb C)$, that is $[f,g]=fg'-f'g$.
\end{proposition}
\begin{proof}
For the sake of convenience of notations, we will equate formal series $a_1z+a_2z^2+...\in\mathbb X_z^0(\mathbb C)$ with their sequence representations $(a_1,a_2,...)\in \mathbb C^{\mathbb N}$ and vice versa. With that convention in mind, one can easily check that for every $f\in\mathbb X_z^0(\mathbb C)$, the tangent space $T_{f}\mathbb X^0_z(\mathbb C)=\mathbb C^{\mathbb N}\equiv \mathbb X^0(\mathbb C)$. The mapping $T\circ$ tangent to $\circ$
is given by (cf. proof of Theorem \ref{smooth}) $T\circ\,((v,h),(w,k))=(v\circ w,(v'\circ w)k+h\circ w)$. Therefore, using the notations from \cite{monastir}, the left-invariant vector fields on $\mathbb X_z^0(\mathbb C)$ can be written as $X_f:\mathbb X_z^0(\mathbb C)\ni g\mapsto g.(z,f):=(g,g'f)$, where $X_f(z):=(z,f)$, and the Lie bracket of two vector fields $X_i\equiv X_{f_i}:g\mapsto (g,\tilde{X_{f_i}}(g):=g'f_i)$ ($i=1,2$) can be defined as $[X_1,X_2]:g\mapsto(g,d_{\tilde{X_1}(g)}\tilde{X_2}(g)-d_{\tilde{X_2}(g)}\tilde{X_1}(g))$. See that, for $(i,j)\in\left\{(1,2),(2,1)\right\}$,
\[
d_{\tilde{X_i}(g)}\tilde{X_j}(g)=\lim\limits_{t\rightarrow 0}\frac{\tilde{X_j}(g+t\tilde{X_i}(g))-\tilde{X_j}(g)}{t}=\lim\limits_{t\rightarrow 0}\frac{(g'+t(g'f_i)')f_j-g'f_j}{t}=g"f_jf_i-g'f_i'f_j,
\]
so $[X_1,X_2](z)=f_1f_2'-f_1'f_2$, which completes the proof.
\end{proof}

Before me move on to some further conslusions, let us introduce a simple
\begin{proposition}
Let $g\in\mathbb X_z^0(\mathbb C)$ and let $S_g$ be a "formal similarity transformation", that is $S_g:\mathbb X_1^{0}(\mathbb C)\ni f\to g\circ f\circ g^{[-1]}\in\mathbb X_1^0(\mathbb C)$. Then $S_g$ is a smooth automorphism of the Fr\'{e}chet-Lie group $\mathbb X_1^{0}(\mathbb C)$. In particular, $\mathbb X_1^0(\mathbb C)$ is a normal subgroup of $\mathbb X_z^0(\mathbb C)$.
\end{proposition}
\begin{proof}
It is easy to check that the operator $S_g$ is bijective and does not change the coefficient $[z^1]f$ of a formal series $f$. Also, for every $f_1,f_2\in\mathbb X_1^0(\mathbb C)$, $S_g(f_1\circ f_2)=g\circ f_1\circ g^{[-1]}\circ g\circ f_2\circ g^{[-1]}=S_g(f_1)\circ S_g(f_2)$ and therefore $S_g(f_1^{[-1]})=(S_g(f_1))^{[-1]}$. The smoothness of $S_g$ is obvious.
\end{proof}
A partial classification of other normal subgroups of $\mathbb X_z^0(\mathbb C)$ is given in the following
\begin{proposition}
Let $G$ be a subgroup of $\mathbb X_z^0(\mathbb C)$ containing no series $f$ such that $[z^1]f$ is a root of unity different from $1$ and containing at least one series $f$ with $[z^1]f\neq 1$. Then $G$ is normal, if and only if there exists such a subgroup $G'$ of $(\mathbb C,\cdot)$ that $G=\left\{a_1z+a_2z^2+...:a_1\in G'\right\}$.
\end{proposition}
\begin{proof}
''$\Rightarrow$'': let $f\in G$, $\alpha:=[z^1]f\neq 1$. Then (\cite{g21}, Prop. 8.1.7.) $\alpha z\in G$, so every series $g$ with $[z^1]g=\alpha$ is in $G$. Also, $\frac{1}{\alpha}z\in G$, so one can conclude that $\mathbb X_1^0\subset G$. Therefore we can write $G=\left\{a_1z+a_2z^2+...:a_1\in A\right\}$, where $A\subset \mathbb C$ (whether a series $g\in G$ depends only on $[z^1]g$). Now, $A$ is a subgroup of $(\mathbb C,\cdot)$ -- indeed, $1\in A$; if $q_1,q_2\in A$, then there exist some formal power series $Q_{1,2}=q_{1,2}z+...\in G$ -- but then $Q_1\circ Q_2\in G$ so $q_1q_2\in A$ and if $q\in A$ then $qz\in G$, so $\frac{1}{q}z\in G$, so $\frac{1}{q}\in A$, which completes this part of the proof. The implication ''$\Leftarrow$'' is obvious.
\end{proof}
The question of how the remaining normal subgroups of $\mathbb X_z^0(\mathbb C)$ -- those containing formal power series $f$ such that $[z^1]f$ is a root of unity -- can be fully described, remains open.
\begin{remark}
One should note that $\mathbb X_1^0(\mathbb C)$ is also not a simple group -- for example, let $n\in\mathbb N$, $n>2$ and let $G=\left\{z+a_nz^n+...:a_n\neq 0\right\}\cup\left\{z\right\}$. It is easy to check, using Lemma \ref{rek}, that $G$ is a subgroup of $\mathbb X_1^0(\mathbb C)$. Now, every $f\in G$ different from $z$ is conjugate to $z+z^n+cz^{2n-1}$ for some $c\in\mathbb C$ (see e.g. \cite{g21}, Prop. 8.1.10.) (that is there exists such $g$ that $f=S_g(z+z^n+cz^{2n-1})$) and therefore for every $g\in\mathbb X_1^0(\mathbb C)$, $S_g(f)$ is conjugate to $z+z^n+cz^{2n-1}$. The claim is now obvious (\cite{g21}, Prop. 8.1.10. and 8.1.11).
\end{remark}
\begin{remark}
See that the Fr\'{e}chet-Lie groups $\mathbb X_z^0(\mathbb C)$, $\mathbb X_1^0(\mathbb C)$ can be treated as formal analogies of the finite-dimensional Lie-groups $GL(n,\mathbb C)$, $SL(n,\mathbb C)$, respectively, where the first coefficients ("$a_1$") of formal power series correspond to the determinants and traces of $n\times n$ matrices. Indeed:
\begin{itemize}
\item by Lemma \ref{rek}, $f\in\mathbb X^0(\mathbb C)$ possesses a composition inverse, if and only if $a_1\neq 0$ ($a_1$ is the formal analogy of the matrix determinant in invertibility conditions);
\item $a_1$ is an invariant of the formal similarity transformation, which is a smooth automorphism of $\mathbb X^0_z(\mathbb C)$ ($\mathbb X^0_1(\mathbb C)$);
\item the Lie group structure of $\mathbb X_z^0(\mathbb C)$ ($\mathbb X_1^0(\mathbb C)$) is analogous to that of $GL(n,\mathbb C)$ ($SL(n\mathbb C)$) in the sense that these sets -- and their corresponding Lie algebras -- are defined by imposing analogous conditions on the respective similarity transformation invariants ($a_1$, matrix determinant and trace); moreover, both $\mathbb X^0_1(\mathbb C)$ and $SL(n,\mathbb C)$ are simply connected, while $\mathbb X^0_z(\mathbb C)$ and $GL(n,\mathbb C)$ can be divided into two connected components -- the set of nonunit power series with $a_1>0$ (the set of $n$-dimensional matrices with positive determinant) and the set of nonunit power series with $a_1<0$ (the set of $n$-dimensional matrices with negative determinant).
\end{itemize}
\end{remark}

\section*{Acknowledgements}
The author is grateful to Professor Dariusz Bugajewski for valuable discussions and a critical reading of the manuscript, and to Professor Piotr Ma\'{c}kowiak for his comments about Theorem \ref{inv}.
\section*{Statements and declarations}
The author declares no competing interests. No new data has been created or analyzed in the presented study.

\end{document}